\documentclass[a4paper,11pt,leqno]{amsart}
\usepackage{amsmath}
\usepackage{amssymb}
\usepackage{amsthm}
\usepackage{indentfirst}
\usepackage[english]{babel}
\usepackage[T1]{fontenc}  
\usepackage{breqn}
\usepackage{cite}    
\usepackage{mathrsfs}                    
\usepackage[colorlinks=true]{hyperref}

\theoremstyle{plain} 
\newtheorem{thm}{Theorem}[section]

\newtheorem{rmk}[thm]{Remark}

\theoremstyle{definition}

\numberwithin{equation}{section}

\author[D.~Mastrostefano]{Daniele Mastrostefano}
\address{University of Warwick, Mathematics Institute, Zeeman Building, Coventry, CV4 7AL, England}
\email{Daniele.Mastrostefano@warwick.ac.uk}


\begin{document}
\title[Positive proportion of short intervals with a given number of primes]
      {Positive proportion of short intervals containing \\
      a prescribed number of primes}\thanks{The author is funded by a Departmental Award and by an EPSRC Doctoral Training Partnership Award.}
      
\begin{abstract}
We will prove that for every $m\geq 0$ there exists an $\varepsilon=\varepsilon(m)>0$ such that if $0<\lambda<\varepsilon$ and $x$ is sufficiently large in terms of $m$ and $\lambda$, then 
$$|\lbrace n\leq x: |[n,n+\lambda\log n]\cap \mathbb{P}|=m\rbrace|\gg_{m,\lambda} x.$$
The value of $\varepsilon(m)$ and the implicit constant on $\lambda$ and $m$ may be made explicit. This is an improvement of an author's previous result. Moreover, we will show that a careful investigation of the proof, apart from some slight changes, can lead to analogous estimates when considering the parameters $m$ and $\lambda$ to vary as functions of $x$ or restricting the primes to belong to specific subsets.
\end{abstract}

\maketitle

\section{Introduction}
\label{sec:1}
Let $\mathbb{P}$ denote the set of prime numbers and fix $\lambda>0$ a real number and $m$ a non-negative integer. The author \cite{M} has recently proved that the proportion of short intervals of the form $[n,n+\lambda\log n]$, for $n\leq x$, containing exactly $m$ primes can be lower bounded by $1/\log x$ if we choose $\lambda$ sufficiently small. More precisely, it was shown that\footnote{Note that the dependence on $\lambda$ of the implicit constant and of $x$ were not stated explicitely in \cite{M}.}
\begin{equation}
\label{eq:1.1}
d_{\lambda,m}(x):=\frac{|\lbrace n\leq x: |[n,n+\lambda\log n]\cap \mathbb{P}|=m\rbrace|}{x}\gg_{m,\lambda} \frac{1}{\log x},
\end{equation}
whenever $0<\lambda<\varepsilon$ for a certain $\varepsilon=\varepsilon(m)>0$ and $x$ large enough in terms of $m$ and $\lambda$.
Under these circumstances, it constitutes a considerable improvement of a previous result of Freiberg\cite{F}, who gave the lower bound $d_{\lambda,m}(x)\gg x^{-\varepsilon'(x)}$, with $\varepsilon'(x)=(\log\log\log\log x)^{2} /(\log\log\log x)$, true for any choice of parameters $\lambda$ and $m$.

The idea behind both those results is that we can make use of the Maynard's sieve method \cite{M1} to find clusters of consecutive primes inside particular sets and then construct short intervals of specific form around them.
Indeed, the Maynard's sieve method allows us to show that any subset of the primes, which is well distributed in arithmetic progressions, contains many elements that are close together. The work of Freiberg showed that the subset of primes which belongs to the image of certain admissible sets of linear functions is well distributed in arithmetic progressions and is suitable for the application of Maynard's results. A combinatorial process is hence used to detect a fixed number among them that are contained in our selected set of intervals.

The major difference between the work of Freiberg and that one of the author is in the way the needed admissible set of linear forms is generated. In the former case, an Erd\H{o}s--Rankin type construction \cite[Lemma 3.3]{F} was considered, which allows us to lower bound the density related to each choice of $\lambda$ and $m$. However, this freedom inevitably forces us to lose precision and obtain weak estimates. In the latter case, the set of linear forms was chosen implicitly by means of the Maynard's sieve, producing in this way better information on the density, only for very small values of $\lambda$.

The aim of the present note is to improve the author's previous work, showing that a better exploration of the last aforementioned approach leads us to generate a positive proportion of short intervals containing a prescribed number of primes. The key idea is that at the start of the process we need to select clusters of primes in which the elements are also well-spaced.

From now on, we will indicate with $m$ a non-negative integer, with $k$ the value $k=C\exp(49m/C')$, for certain suitable constants $C,C'>0$, and with $\lambda$ a positive real number smaller than $\varepsilon=\varepsilon(k):=k^{-4}(\log k)^{-2}$. The result is the following.
\begin{thm}
\label{thm 1.1}
We have
\begin{equation}
\label{eq:1.2}
d_{\lambda,m}(x)\gg \lambda^{k+1}e^{-Dk^{4}\log k},
\end{equation}
for a certain absolute constant $D>0$, if $x$ is sufficiently large in terms of $m$ and $\lambda$.
\end{thm}

It is interesting to note that, from a heuristic point of view, we expect a positive proportion result for all the short intervals of the form $[n,n+\lambda\log n]$ (and for all the non-negative integers $m$). More precisely, we conjecture that 
$$d_{\lambda,m}(x)\backsim \frac{\lambda^{m} e^{-\lambda}}{m!},\ \textrm{as}\ x\rightarrow\infty$$
for every $\lambda$ and $m$ (see for instance the expository article \cite{S} of Soundararajan for further discussions). 

The strength of the Maynard's sieve is the flexibility, that makes it applicable to counting primes in sparser subsets as well. In fact, the same proof that leads to Theorem \ref{thm 1.1} can be overall adapted to study a variety of different situations, in which for instance we restrict the primes to lye on an arithmetic progression or allow for uniformity of the parameters $\lambda$ and $m$. The results are as following.
\begin{thm}
\label{thm 1.2}
Let $x$ be sufficiently large in terms of $m$ and $\lambda$. Suppose that $q\leq f(x)$ is a positive integer, with $(\log x)/f(x)\rightarrow\infty$, as $x\rightarrow\infty$. Take $0\leq a<q$ with $(a,q)=1$. Then, we have
\begin{equation}
\label{eq:1.3}
d_{\lambda,m}^{a,q}(x)\gg \frac{\lambda^{k+1}e^{-Dk^{4}\log k}}{q^{k+1}},
\end{equation}
for a certain $D>0$, where $d_{\lambda,m}^{a,q}(x)$ is defined as in \eqref{eq:1.1} but with $\mathbb{P}$ replaced by $\mathbb{P}_{a,q}$, being the intersection of $\mathbb{P}$ with $a\pmod{q}$.
\end{thm}
\begin{thm}
\label{thm 1.3}
Fix $\epsilon_1>0$ a small parameter and $0<\epsilon_{2}<1$. Let $x\geq x_{0}(\epsilon_1,\epsilon_2)$, $m\leq \epsilon_1\log\log x$ and $\lambda\geq (\log x)^{\epsilon_2-1}$, obeying to the relations $k^4(\log k)^2 \lambda\leq 1$ and $\lambda>k\log k (\log x)^{-1}$. Then, the estimate \eqref{eq:1.2} continues to hold.
\end{thm}
\begin{thm}
\label{thm 1.4}
Let $\mathbb{K}/\mathbb{Q}$ be a Galois extension of $\mathbb{Q}$ with discriminant $\Delta_{\mathbb{K}}$. There exist constants $C_{\mathbb{K}},C_{\mathbb{K}}'>0$ depending only on $\mathbb{K}$ such that the following holds. Let $\mathcal{C} \subset Gal(\mathbb{K}/\mathbb{Q})$ be a conjugacy class in the Galois group of $\mathbb{K}/\mathbb{Q}$, and let 
$$\mathcal{P} = \lbrace p\ \textrm{prime} : p\nmid \Delta_{\mathbb{K}}, \left[\frac{\mathbb{K}/\mathbb{Q}}{p}\right]=\mathcal{C}\rbrace,$$
where $\left[\frac{\mathbb{K}/\mathbb{Q}}{.}\right]$ denotes the Artin symbol. Let $m\in\mathbb{N}$, $k=C_{\mathbb{K}}'\exp(C_{\mathbb{K}} m)$ and $\lambda<\varepsilon$. Then, we have
\begin{equation}
\label{eq:1.4}
d_{\lambda,m}^{\mathbb{K}}(x)\gg\lambda^{k+1}e^{-Dk^{4}\log k},
\end{equation}
provided $x \geq x_0(\mathbb{K},\lambda,m)$, where $d_{\lambda,m}^{\mathbb{K}}(x)$ is defined as in \eqref{eq:1.1} except that $\mathbb{P}$ is replaced by $\mathcal{P}$.
\end{thm}
If we consider values of $\lambda$ slightly bigger than $k^{-4}(\log k)^{-2}$, a little variation of the sieve method used to prove Theorem \ref{thm 1.1} leads to the following improvement on the Freiberg bound in \cite{F}.
\begin{thm}
\label{thm 1.5}
For every non-negative integer $m$ and positive real number $\lambda$ smaller than $k^{-1}(\log k)^{-1}$, with $k$ the value $k=C\exp(49m/C')$, for suitable constants $C,C'>0$, we have 
\begin{equation}
\label{eq:1.5}
d_{\lambda,m}(x)\gg \frac{\lambda e^{-Dk^{4}\log k}}{(\log x)^{k}},
\end{equation}
for a certain $D>0$, if $x$ is sufficiently large in terms of $m$ and $\lambda$.
\end{thm}
\section{Notations}
Throughout, $\mathbb{P}$ denotes the set of all primes, $\mathbf{1}_{S}: \mathbb{N}\rightarrow \lbrace 0,1\rbrace$ the indicator function of a set $S\subset \mathbb{N}$ and $p$ a prime. As usual, $\varphi$ will denote the Euler totient function and $(m,n)$ the greatest common divisor of integers $n$ and $m$. We will always denote with $x$ a sufficiently large real number. By $o(1)$ we mean a quantity that tends to $0$ as $x$ tends to infinity. The expressions $A=O(B), A\ll B, B\gg A$ denote that $|A|\leq c|B|$, where $c$ is some positive (absolute, unless stated otherwise) constant.

In the following we will always consider admissible $k$-tuples of linear forms $\lbrace gn+h_{1},...,gn+h_{k}\rbrace$, where $0\leq h_{1}<h_{2}<...<h_{k}<\lambda\log x$, $k$ a sufficiently large integer and $g$ a positive integer, coprime with $B$, squarefree and such that $\log x<g\leq 2\log x$. Here, $B=1$ or $B$ is a prime with $\log\log x^{\eta}\ll B\ll x^{2\eta},$ where we put $\eta:=c/500k^{2}$ with $0<c<1$. As usual, a finite set $\mathcal{L}:=\lbrace L_{1},...,L_{k}\rbrace$ of linear functions is admissible if the set of solutions modulo $p$ to $L_ {1}(n)\cdots L_{k}(n)\equiv 0\pmod{p}$ does not form a complete residue system modulo $p$, for any prime $p$. In our case, in which $L_{i}(n)=gn+h_{i}$, for every $i=1,...,k$, we may infer that the set $\lbrace L_{1},...,L_{k}\rbrace$ is admissible if and only if the set $\mathcal{H}:=\lbrace h_{1},...,h_{k}\rbrace$ it is, in the sense that the elements $h_{1},...,h_{k}$ do not cover all the residue classes modulo $p$, for any prime $p$.

The proof of Theorem \ref{thm 1.1} (and of its variations \ref{thm 1.2}--\ref{thm 1.5}) follows by mimicking that one in \cite{M}, taking into account a new crucial assumption on the set of linear forms we will work with. We will briefly rewrite the main estimates and passages already containted in the proof in \cite{M}, highlighting the main differences and the new computations. In particular, several notations will not be introduced here because not essential for the general understanding of the argument or already present in \cite{M}.
\section{Application of the Maynard's Sieve}
As at the start of \cite[Section 3]{M}, and following the notations there introduced, we define the double sum
\begin{equation}
\label{eq:3.1}
S=\sum^{*}_{\mathcal{H}}\sum_{x<n\leq 2x}S(\mathcal{H},n),
\end{equation}
where 
\begin{equation}
\label{eq:3.2}
S(\mathcal{H},n)=\bigg( \sum_{i=1}^{k} \textbf{1}_{\mathbb{P}}(gn+h_{i})-m-k\sum_{i=1}^{k}\sum_{\substack{p|gn+h_{i}\\ p\leq x^{\rho}, p\nmid B}}1 -k\sum_{\substack{ h\leq 5\lambda\log x\\ (h,g)=1\\h\not\in\mathcal{H}}} \textbf{1}_{S(\rho, B)}(gn+h)\bigg)w_{n}(\mathcal{H}).
\end{equation}
Here $\Sigma^{*}_{\mathcal{H}}$ means that the sum is over all the admissible sets $\mathcal{H}$ such that $0\leq h_{1}<h_{2}<...<h_{k}<\lambda\log x$ and $|h_i-h_j|>\frac{\lambda\log x}{C_{0}}$, for any $1\leq i\neq j\leq k$, where $C_0$ will be chosen later. Note that, unlike in \cite{M}, in the innermost sum in \eqref{eq:3.1} we now have $m$ instead of $m-1$. 

Following closely the discussion at the beginning of \cite[Section 3]{M}, we deduce that
\begin{equation}
\label{eq:3.3}
S\ll k(\log x)^{2k} \exp (O(k/\rho))|I(x)|,
\end{equation}
where now the set $I(x)$ contains intervals of the form  $[gn,gn+5\lambda\log x]$, for $x<n\leq 2x$, with  the property that $|[gn,gn+5\lambda\log x]\cap \mathbb{P}|=|\lbrace gn+h_{1},...,gn+h_{k}\rbrace\cap\mathbb{P}|\geq m+1$, for a unique admissible set $\mathcal{H}$ such that $0\leq h_{1}<...<h_{k}<\lambda\log x$ and $|h_i-h_j|>\frac{\lambda\log x}{C_{0}}$, for any $1\leq i\neq j\leq k$. We recall that the intervals in $I(x)$ are pairwise disjoint, if for instance $\lambda<1/5$.\\
We need also a lower bound for $S$. Using \cite[Proposition 2.1, 2.2--2.5]{M}, we find
\begin{equation}
\label{eq:3.4}
S\geq\sum_{\mathcal{H}}^{*}\bigg[ (1+o(1))\frac{B^{k-1}}{\varphi(B)^{k-1}•}\mathfrak{S}_{B}(\mathcal{H})(\log R)^{k+1} J_{k}\sum_{i=1}^{k}\frac{\varphi(g)}{g•}\sum_{x<n\leq 2x}\textbf{1}_{\mathbb{P}}(gn+h_{i})
\end{equation}
$$-m(1+o(1))\frac{B^{k}}{\varphi(B)^{k}•}\mathfrak{S}_{B}(\mathcal{H})x(\log R)^{k} I_{k}+O\left(\rho^{2}k^{6}(\log k)^{2}\frac{B^{k}}{\varphi(B)^{k}•}\mathfrak{S}_{B}(\mathcal{H})x(\log R)^{k} I_{k}\right)$$
$$+O\left(k\frac{B^{k}}{\varphi(B)^{k}•}\mathfrak{S}_{B}(\mathcal{H})x(\log R)^{k-1} I_{k}\right)+O\bigg(\frac{k}{\rho}\frac{B^{k+1}}{\varphi(B)^{k+1}•}\mathfrak{S}_{B}(\mathcal{H})x(\log R)^{k-1} I_{k}\sum_{\substack{h\leq 5\lambda\log x\\(h,g)=1\\ h\not\in \mathcal{H}}}\frac{\Delta_{\mathcal{L}}}{\varphi(\Delta_{\mathcal{L}•})}\bigg)\bigg].$$
By the inequality (2.7) in \cite{M}, we have
$$\frac{\varphi(B)}{B•}\frac{\varphi(g)}{g•}\sum_{i=1}^{k}\sum_{x<n\leq 2x}\textbf{1}_{\mathbb{P}}(gn+h_{i})>\frac{kx}{2\log x•}.$$
Note that the hypotheses of \cite[Theorem 2.2]{M} are satisfied. Using this estimate together with \cite[Lemma 3.1]{M}, and choosing $\rho:=k^{-3}(\log k)^{-1}$, we find
\begin{equation}
\label{eq:3.5}
S\geq\sum_{\mathcal{H}}^{*}\frac{B^{k}}{\varphi(B)^{k}}\mathfrak{S}_{B}(\mathcal{H})x(\log R)^{k}\bigg[(1+o(1))kJ_{k}\frac{\log R}{2\log x}-mI_{k}(1+o(1))+O(I_{k})
\end{equation}
$$+O(k I_{k}(\log R)^{-1})+O\left(k^{4}(\log k) I_{k}(\log R)^{-1}\lambda\log x(\log k) \right)\bigg].$$
We remark here that the aforementioned Theorem 2.2 and Lemma 3.1 need $x$ to be large enough with respect to $k$ and $\lambda$.
Now, by \cite[Proposition 2.1, 2.6]{M} we know that $J_{k}\geq C' \frac{\log k}{k}I_{k},$ for a certain $C'>0$. We should consider $k$ sufficiently large in terms of $m$. For example, we may take $k:=C\exp(49m/C')$, with $C>0$. Choosing $\lambda\leq\varepsilon$, with $\varepsilon=\varepsilon(k):=k^{-4}(\log k)^{-2}$, and taking $x$ and $C$ suitably large, we may conclude that
\begin{equation}
\label{eq:3.6}
S\gg\sum_{\mathcal{H}}^{*}\frac{B^{k}}{\varphi(B)^{k}}\mathfrak{S}_{B}(\mathcal{H})x(\log R)^{k}I_{k}.
\end{equation}
By the estimates \cite[Proposition 2.1, 2.6]{M} we know that $I_{k}\gg (2k\log k)^{-k}$ and $\mathfrak{S}_{B}(\mathcal{H})\gg \exp(-C_{1}k)$, for a certain $C_{1}>0$. Remember also that $R=x^{1/24}$. Finally, we may certainly use $\frac{B^{k}}{\varphi(B)^{k}}\geq 1$. Inserting all of these in \eqref{eq:3.6}, we obtain
\begin{equation}
\label{eq:3.7}
S\gg x(\log x)^{k}e^{-C_{2}k^{2}}\sum_{\mathcal{H}}^{*} 1,
\end{equation}
for a suitable constant $C_{2}>0$. Thus, we are left with obtaining a lower bound for the sum in \eqref{eq:3.7}.
We greedily sieve the interval $[0,\lambda\log x]$, by removing for each prime $p\leq k$ in turn any elements from the residue class modulo $p$ which contains the fewest elements. The resulting set $\mathcal{A}$, say, has size
$$|\mathcal{A}|\geq \lambda\log x\prod_{p\leq k}\left(1-\frac{1}{p}\right)\geq c'\frac{\lambda\log x}{(\log k)},$$
by Mertens's theorem, with $c'>0.$ 

Any choice of $k$ distinct $h_{i}$ from $\mathcal{A}$ will constitute an admissible set $\mathcal{H}=\lbrace h_{1},...,h_{k}\rbrace$ such that $0\leq h_{1}<h_{2}<...<h_{k}< \lambda\log x$. Now, we count how many of them have the property that $|h_i-h_j|>\frac{\lambda\log x}{C_{0}}$, for any $1\leq i\neq j\leq k$. Certainly, we can choose $h_1$ in $|\mathcal{A}|$-ways. Let us call $\mathcal{A}_1:=\mathcal{A}$ and define 
$$\mathcal{A}_2:=\mathcal{A}_1\setminus \mathcal{A}_1\cap [h_1-\lfloor \lambda\log x/C_{0}\rfloor, h_1+\lfloor \lambda\log x/C_{0} \rfloor].$$
We will pick $h_2\in\mathcal{A}_2$, having then $|\mathcal{A}_2|\geq |\mathcal{A}_1|-2\lfloor \lambda\log x/C_{0}\rfloor$ possibilities.
Iterating this process, we can count the number of admissible choices for any $h_i$ until $h_k$, which will be an element in 
$$\mathcal{A}_k:=\mathcal{A}_{k-1}\setminus \mathcal{A}_{k-1}\cap [h_{k-1}-\lfloor \lambda\log x/C_{0}\rfloor, h_{k-1}+\lfloor \lambda\log x/C_{0} \rfloor],$$
which will have cardinality $|\mathcal{A}_k|\geq |\mathcal{A}_1|-2(k-1)\lfloor \lambda\log x/C_{0}\rfloor$.

In conclusion, for our particular choice of admissible sets we have at least a number of possibilities equals to
\begin{equation}
\label{eq:3.8}
\frac1{k!}\prod_{i=1}^{k}|\mathcal{A}_i|\geq \frac1{k^k}\prod_{i=1}^{k}\left(|\mathcal{A}_1|-2(i-1)\lfloor \lambda\log x/C_{0}\rfloor\right)\geq \frac1{k^{k}}\left( c'\frac{\lambda\log x}{(\log k)}-2k\frac{\lambda\log x}{C_{0}}\right)^{k}.
\end{equation}
Let's take $C_0=C_0(k):=4k(\log k)/c'$. We immediately see that \eqref{eq:3.8} becomes $\gg \lambda^{k}e^{-C_3 k^2}(\log x)^{k},$ which leads to $S\gg \lambda^{k} e^{-C_{4} k^{2}}x(\log x)^{2k},$ for certain constants $C_3, C_{4}>0$. Finally, by combining \eqref{eq:3.3} with the above information on $S$ we obtain
\begin{equation}
\label{eq:3.9}
|I(x)|\gg \lambda^{k}e^{-C_{5}k^{4}\log k}x,
\end{equation}
with an absolute constant $C_{5}>0$.
\section{Modification of the combinatorial process}
Consider an interval $I\in I(x)$. There exist an integer $x<n\leq 2x$ and an admissible set $\mathcal{H}$, with $0\leq h_{1}<h_{2}<...<h_{k}< \lambda\log x$ and $|h_i-h_j|>\lambda\log x/C_0,$ for any $1\leq i\neq j\leq k$, such that $I=[gn, gn+5\lambda\log x]$ and 
$$|[gn,gn+5\lambda\log x]\cap \mathbb{P}|=|\lbrace gn+h_{1},...,gn+h_{k}\rbrace\cap\mathbb{P}|\geq m+1.$$
In order to avoid having a trivial gap between the elements of $\mathcal{H}$ we ask for $x$ to be sufficiently large with respect to $\lambda$ and $k$. Let us define 
\begin{equation}
\label{eq:4.1}
I_{j}=[N_{j}, N_{j}+\lambda \log N_{j}],\ \ N_{j}=gn+j,
\end{equation}
for $j=0,...,\lfloor \lambda\log N_{0}\rfloor$. We recall here the following properties of the intervals $I_j$, that are stated and proved in details in \cite{M}:\\
\\
\begin{tabular}{@{}l}
\textbf{1)} for any such $j$ we have $I_{j}\subseteq I$;\\
\textbf{2)} for the choice $j=h_{1}$ we find that $I_{j}\cap \lbrace gn+h_{1},...,gn+h_{k}\rbrace=\lbrace gn+h_{1},...,gn+h_{k}\rbrace$;\\
\textbf{3)} for the value $j=\lfloor \lambda\log N_{0}\rfloor$ we have $I_{j}\cap\lbrace gn+h_{1},...,gn+h_{k}\rbrace=\emptyset$;\\
\textbf{4)} if $|I_{j}\cap \mathbb{P}|<|I_{j+1}\cap\mathbb{P}|$, for a certain $j$, then $|I_{j+1}\cap\mathbb{P}|=|I_{j}\cap\mathbb{P}|+1.$\\
\end{tabular}\
\\
\\
Now, let's define 
$$\tilde{j}:=\max\lbrace 0\leq j\leq \lfloor \lambda\log N_{0}\rfloor: |I_j\cap \mathbb{P}|\geq m+1\rbrace.$$
Note that we necessarily have $N_{\tilde{j}}=gn+\tilde{j}$ being prime. Consequently, this implies $|I_{\tilde{j}+1}\cap \mathbb{P}|=m$, but from our assumption on $\mathcal{H}$ it actually derives that 
$$|I_{\tilde{j}+l}\cap \mathbb{P}|=m,\ \textrm{for any}\ 1\leq l\leq \lfloor \lambda\log x/C_0\rfloor.$$
This is equivalent to say that we have found $\lfloor \lambda\log x/C_0\rfloor$-different intervals $[N, N+\lambda\log N]$ containing exactly $m$ primes, with $N<5x\log x$, if $x$ is sufficiently large. Together with the lower bound \eqref{eq:3.9}, we have obtained that for every $m\geq 0$ and for each $\lambda\leq\varepsilon$,
\begin{equation}
\label{eq:4.2}
|\lbrace N\leq 5x\log x: |[N, N+\lambda\log N]\cap\mathbb{P}|=m\rbrace|\gg \frac{\lambda^{k+1}}{k\log k}e^{-C_{5}k^{4}\log k}x\log x,
\end{equation}
which is equivalent to
\begin{equation}
\label{eq:4.3}
|\lbrace N\leq X: |[N, N+\lambda\log N]\cap\mathbb{P}|=m\rbrace|\gg \lambda^{k+1}e^{-C_{6}k^{4}\log k}X,
\end{equation}
when $X$ is large enough in terms of $\lambda$ and $k$, for a certain constant $C_6>0$, which proves Theorem \ref{thm 1.1}. 
\section{Concluding Remarks}
\subsection{Explicit constants}
Since we let $k=C\exp(49m/C')$, with $C,C'$ as above, note that we can rewrite the final estimate using only the dependence in $\lambda$ and $m$.\\
Remembering the choice of $\varepsilon(m)$, it is immediate to see the following interdependence between $\lambda$ and $m$, given by 
\begin{equation}
\label{eq:5.1}
\lambda(49m+c_1)^{2}\exp(196c_2 m)\ll 1
\end{equation}
for certain constants $c_1,c_2>0.$
\begin{rmk}
Notice that, if we might take $k=m+1$ and we were able to improve the constants in $k$ in the sieve method, then we would end up with an explicit constant in \eqref{eq:4.3} that almost matches the expected one, for values of $\lambda$ close to $0$.
\end{rmk}
\subsection{The case of primes in arithmetic progressions}
Suppose that $q$ is a squarefree positive integer, coprime with $B$ and $q\leq f(x)$ with $(\log x)/f(x)\rightarrow\infty$, as $x\rightarrow\infty$. Take $0\leq a<q$ with $(a,q)=1$. In order to extend the result of Theorem \ref{thm 1.1} to this situation we go over again its proof. In particular, in \eqref{eq:3.1} we average now over admissible sets $\mathcal{H}=\lbrace h_{1},...,h_{k}\rbrace$ such that $0\leq h_{1}:=a+qb_1<h_{2}:=a+qb_2<...<h_{k}:=a+qb_k<\lambda\log x$, and $|b_i-b_j|>\frac{\lambda\log x}{qC_{0}}$, for any $1\leq i\neq j\leq k$. Moreover, we need to take $g$ as a squarefree multiple of $q$, coprime with $B$ and such that $\log x<g\leq 2\log x$. Therefore, such set of linear functions satisfies the hypotheses of \cite[Proposition 2.1]{M} and \cite[Theorem 2.2]{M}, and the images of all its elements lye on the arithmetic progression $a\pmod{q}$.\\
In particular, we obtain analogously to what done in section 3 that
$$|I(x)|\gg \frac{\lambda^{k}e^{-C_8 k^{4}\log k}}{q^{k}}x,$$
if $x$ is sufficiently large in terms of $\lambda$ and $k$, for a suitable constant $C_8>0$.\\
Here, using the notation $\mathbb{P}_{a,q}$ to indicate the primes in the arithmetic progression $a\pmod{q}$, the set $I(x)$ contains intervals of the form $[gn,gn+5\lambda\log x]$, for $x<n\leq 2x$, with $g\equiv 0\pmod{q}$ and the property that $|[gn,gn+5\lambda\log x]\cap \mathbb{P}|=|\lbrace gn+h_{1},...,gn+h_{k}\rbrace\cap\mathbb{P}_{a,q}|\geq m+1$, for a unique admissible set $\mathcal{H}=\lbrace h_{1},...,h_{k}\rbrace$ such that $0\leq h_{1}:=a+qb_1<h_{2}:=a+qb_2<...<h_{k}:=a+qb_k<\lambda\log x$ and $|b_i-b_j|>\frac{\lambda\log x}{qC_{0}}$, for any $1\leq i\neq j\leq k$.\\
Following the computations done in section 4, we obtain
$$|\lbrace N\leq X: |[N, N+\lambda\log N]\cap\mathbb{P}_{a,q}|=m\rbrace|\gg \frac{\lambda^{k+1}e^{-C_9 k^{4}\log k}}{q^{k+1}}X,$$
when $X$ is sufficiently large in terms of $\lambda$ and $k$, for a suitable absolute constant $C_9>0$, which proves Theorem \ref{thm 1.2}. Indeed, the restriction on $q$ to be squarefree and coprime with $B$ can be removed at the cost to slightly modify the proof of \cite[Theorem 2.2]{M}. In particular, at the start of its proof we need to change $B$ with the largest prime factor of $\tilde{l}$ coprime with $g$, with $\tilde{l}$ being the modulus of a possible exceptional character among all the primitive Dirichlet characters $\chi\bmod{l}$ to moduli $l\leq x^{2\eta}$.
\subsection{The case of uniform parameters}
In section 3 we applied \cite[Proposition 2.1]{M}, which is a specific case of \cite[Proposition 6.1]{M1}, in which a uniformity in $k\leq (\log x)^{1/5}$, say, is allowed. A careful examination of \cite[Theorem 2.2]{M}, and of the computations done in Section 3 and 4 in the present paper and in \cite{M}, shows that the estimate \eqref{eq:1.2} continues to hold also when $m\leq \epsilon_1\log\log x$ and $\lambda\geq (\log x)^{\epsilon_2-1}$ satisfy \eqref{eq:5.1} together with $\lambda>k\log k(\log x)^{-1}$. Here, $\epsilon_1$ is a fixed sufficiently small constant (e.g. smaller than $C'/294$) and $0<\epsilon_2<1$.
\begin{rmk}
Notice that, in the case in which $\lambda$ goes to $0$ together with $x$, and $m$ and $\lambda$ vary in the range defined above, the Cram\'er model used in \cite{S} still gives us an expected asymptotic value for $d_{\lambda,m}$, which now takes the form 
$$d_{\lambda,m}(x)\sim \frac{\lambda^{m}}{m!},\ \textrm{as}\ x\rightarrow\infty.$$
Obviously, since the constant in $m$ in the lower bound \eqref{eq:1.2} is not optimal, the result of Theorem \ref{eq:1.3} now will be far away from what the model suggests.
\end{rmk}
\subsection{The case of primes in Chebotarev sets}
As already mentioned, we have only used so far a very special case of \cite[Proposition 6.1]{M1}. In particular, it is meaningful to observe that we can replace the set of all the primes with a fairly smaller one, as long as it verifies a suitable variant of \cite[Theorem 2.2]{M}. More specifically, we would like to concentrate on the so called primes in Chebotarev sets.\\
Let $\mathbb{K}/\mathbb{Q}$ be a Galois extension of $\mathbb{Q}$ with discriminant $\Delta_{\mathbb{K}}$. Let $\mathcal{C} \subset Gal(\mathbb{K}/\mathbb{Q})$ be a conjugacy class in the Galois group of $\mathbb{K}/\mathbb{Q}$, and let 
$$\mathcal{P} = \lbrace p\ \textrm{prime} : p\nmid \Delta_{\mathbb{K}}, \left[\frac{\mathbb{K}/\mathbb{Q}}{p}\right]=\mathcal{C}\rbrace,$$
where $\left[\frac{\mathbb{K}/\mathbb{Q}}{.}\right]$ denotes the Artin symbol. Fix $m\in\mathbb{N}$, $k=C_{\mathbb{K}}' \exp(C_{\mathbb{K}} m)$, for suitable $C_{\mathbb{K}},C_{\mathbb{K}}'>0$, and $\lambda<\varepsilon$. Finally, let $\log x<g\leq 2\log x$ be a squarefree number with $(g,\Delta_{\mathbb{K}})=1$, bearing in mind that now $B=\Delta_{\mathbb{K}}$, and consider admissible sets $\mathcal{H}$ of the usual form.\\
Murty and Murty proved in their main theorem in \cite{MM} that the primes in $\mathcal{P}$ are well distributed among arithmetic progressions of moduli $q\leq x^{\theta}$, with $\theta<\min(1/2,2/|G|)$, and such that $\mathbb{K}\cap \mathbb{Q}(\zeta_{q})=\mathbb{Q}$. An adaptation of the argument present in the proof of \cite[Theorem 2.2]{M} leads to prove the second estimate stated there, where $\mathbb{P}$ is replaced by $\mathcal{P}$ and the sum over $q$ is over all the moduli $q\leq x^{\theta/4}$, satisfying the algebraic condition described above. Regarding the first estimate in \cite[Theorem 2.2]{M}, we have
$$\frac{1}{k}\frac{B}{\varphi(B)}\frac{\varphi(g)}{g}\sum_{i=1}^{k}\sum_{x<n\leq 2x}\textbf{1}_{\mathcal{P}}(gn+h_i)\geq(1 + o(1))\frac{\Delta_{\mathbb{K}}}{\varphi(\Delta_{\mathbb{K}})}\frac{|\mathcal{C}|}{|G|}\frac{x}{\log x},$$
which essentially follows from the Chebotarev density theorem. Working as in section 3, we can find
$$|I(x)|\gg \lambda^{k}e^{-C_{10} k^{4}\log k}x,$$
if $x\geq x_{0}(\mathbb{K},\lambda,m)$, for a suitable constant $C_{10}>0$.\\
Here, the set $I(x)$ contains interval of the form $[gn,gn+5\lambda\log x]$, for $x<n\leq 2x$ and $\log x<g\leq 2\log x$, having the property that $|[gn,gn+5\lambda\log x]\cap \mathbb{P}|=|\lbrace gn+h_{1},...,gn+h_{k}\rbrace\cap\mathcal{P}|\geq m+1$, for a unique admissible set $\mathcal{H}$ such that $0\leq h_{1}<h_{2}<...<h_{k}<\lambda\log x$ and $|h_i-h_j|>\frac{\lambda\log x}{C_{0}}$, for any $1\leq i\neq j\leq k$.\\
Following the computations done in section 4, we obtain
$$|\lbrace N\leq X: |[N, N+\lambda\log N]\cap\mathcal{P}|=m\rbrace|\gg \lambda^{k+1}e^{-C_{11} k^{4}\log k}X,$$
when $X\geq X_{0}(\mathbb{K},\lambda,m)$, for a suitable absolute constant $C_{11}>0$, which proves Theorem \ref{thm 1.4}.
\subsection{The case of slightly bigger values of $\lambda$}
Let us fix an admissible $k$-tuple of linear functions $\mathcal{L}=\lbrace gn+h_{1},...,gn+h_{k}\rbrace$ with the usual form. We replace the last sum in parenthesis in \eqref{eq:3.2} with the following one
\begin{equation*}
\sum_{\substack{ h\leq 5\lambda\log x\\ (h,g)=1\\h\not\in\mathcal{H}}} \textbf{1}_{S(1/80, 1)}(gn+h)
\end{equation*} 
and we remove the average over $\mathcal{H}$ in \eqref{eq:3.1}. With these variations in mind, it is immediate to see that \eqref{eq:3.3} still continues to hold, but now we can only say that for every interval $I\in I(x)$ there exists an integer $x<n\leq 2x$ such that $I=[gn,gn+5\lambda\log x]$ and 
$$|[gn,gn+5\lambda\log x]\cap \mathbb{P}|=|\lbrace gn+h_{1},...,gn+h_{k}\rbrace\cap\mathbb{P}|\geq m+1.$$
Arguing as in section 3 with the opportune variations, but essentially carrying over all the computations, we deduce that
\begin{equation}
\label{eq:5.2}
S\gg x(\log x)^{k}e^{-C_{12}k^{2}},\ |I(x)|\gg e^{-C_{13}k^{4}\log k}\frac{x}{(\log x)^{k}},
\end{equation}
for suitable $C_{12},C_{13}>0$. The only key difference in proving \eqref{eq:5.2} is that the last big-O in \eqref{eq:3.4} now assumes the shape
$$O\bigg(80k\frac{B^{k}}{\varphi(B)^{k}}\mathfrak{S}_{B}(\mathcal{H})x(\log R)^{k-1} I_{k}\sum_{\substack{h\leq 5\lambda\log x\\(h,g)=1\\ h\not\in \mathcal{H}}}\frac{\Delta_{\mathcal{L}}}{\varphi(\Delta_{\mathcal{L}•})}\bigg).$$
Consequently, this modifies also the last big-O in \eqref{eq:3.5}, which will be
$$O\left(k(\log k)I_{k}(\log R)^{-1}\lambda\log x \right)=O(I_{k}),\ \textrm{if}\ \lambda<\frac{1}{k\log k}.$$
The rest of the argument goes through as always and we conclude that
\begin{equation*}
|\lbrace N\leq X: |[N, N+\lambda\log N]\cap\mathbb{P}|=m\rbrace|\gg \lambda e^{-C_{14}k^{4}\log k}\frac{X}{(\log X)^{k}},
\end{equation*}
when $X$ is large enough in terms of $\lambda$ and $k$, for a certain $C_{14}>0$, which proves Theorem \ref{thm 1.5}. 
\begin{rmk}
We would like to observe that many of the variables and parameters have not been chosen in the best possible way, since finding their precise range of definition is not in the spirit of the paper and does not considerably improve the final results. We refer to \cite{T} for several arithmetic consequences of finding primes of a given splitting type and note that they may be translated in our context.
Finally, we would like to point out that we are able to mix up the results presented in this section, paying attention to the possible relations between the different parameters.
\end{rmk}

\end{document}